\documentclass[12pt]{amsart} 
\usepackage{amsmath,amssymb,amsthm}
\usepackage[all]{xy}

\newcommand{\Sym}{\mathfrak{S}}
\newcommand{\C}{\mathbb{C}}
\newcommand{\N}{\mathbb{N}}
\newcommand{\Z}{\mathbb{Z}}

\newcommand{\GL}{\mathsf{GL}}
\newcommand{\gl}{\mathfrak{gl}}

\newcommand{\Hom}{\operatorname{Hom}}
\newcommand{\End}{\operatorname{End}}
\newcommand{\sgn}{\operatorname{sgn}}
\newcommand{\ann}{\operatorname{ann}}
\newcommand{\ind}{\operatorname{ind}}

\newcommand{\dom}{\unrhd}
\newcommand{\sdom}{\rhd}
\renewcommand{\le}{\leqslant}
\renewcommand{\ge}{\geqslant}

\newcommand{\cf}{\operatorname{cf}}
\newcommand{\shape}{\operatorname{shape}}
\newcommand{\ch}{\operatorname{ch}}
\newcommand{\transpose}{\mathsf{T}}

\newtheorem{thm}{Theorem}
\newtheorem{lem}[thm]{Lemma}
\newtheorem{prop}[thm]{Proposition}
\newtheorem{cor}[thm]{Corollary}
\theoremstyle{definition}
\newtheorem{rmk}[thm]{Remark}

\allowdisplaybreaks
\numberwithin{equation}{section}

\begin{document}
\title
{Schur--Weyl duality in positive characteristic}
\author{Stephen Doty}
\address{Mathematics and Statistics, 
Loyola University Chicago, 
Chicago, Illinois 60626 U.S.A.}
\email{sdoty@luc.edu} \thanks{The author is grateful to Jun~Hu for
bringing reference \cite{Harterich} to his attention, and to the
referee for useful suggestions.}

\begin{abstract} 
  Complete proofs of Schur--Weyl duality in positive characteristic
  are scarce in the literature. The purpose of this survey is to write
  out the details of such a proof, deriving the result in positive
  characteristic from the classical result in characteristic zero,
  using only known facts from representation theory.
\end{abstract}
\maketitle

\section{Introduction}\noindent

Given a set $A$ write $\Sym_A$ for the symmetric group
on $A$, i.e., the group of bijections of $A$. For $\sigma \in \Sym_A$
and $a \in A$ we always write $a\sigma$ for the image of $a$ under
$\sigma$. In other words, we choose to write maps in $\Sym_A$ on the
right of their argument. This means that $\sigma \tau$ (for $\sigma,
\tau \in \Sym_A$) is defined by $a(\sigma\tau)=(a\sigma)\tau$.

We will write $\Sym_r$ as a shorthand for $\Sym_{\{1,\dots,r\}}$.

Consider the group $\Gamma=\GL(V)$ of linear automorphisms on an
$n$-dimensional vector space $V$ over a field $K$. We write elements
$g \in \Gamma$ on the left of their argument. (Indeed, maps are
generally written on the left in this article, except when they belong
to a symmetric group.) The given action $(g, v) \to g(v)$ of $\Gamma$
on $V$ induces a corresponding action on a tensor power $V^{\otimes
  r}$, with $\Gamma$ acting the same in each tensor position: $g(u_1
\otimes \cdots \otimes u_r) = (g(u_1))\otimes \cdots \otimes
(g(u_r))$, for $g\in \Gamma$, $u_i \in V$.  Evidently the action of
$\Gamma$ commutes with the ``place permutation'' action of $\Sym_r$,
acting on $V^{\otimes r}$ on the right via the rule $(u_1\otimes
\cdots \otimes u_r) \sigma = u_{1\sigma^{-1}} \otimes \cdots \otimes
u_{r\sigma^{-1}}$. In this action, a vector that started in tensor
position $i\sigma^{-1}$ ends up in tensor position $i$, thus a
vector that started in tensor position $i$ ends up in tensor position
$i\sigma$.

We write $KG$ for the group algebra of a group $G$. The fact that the
two actions commute means that the corresponding representations
\begin{equation}\label{eq:repdefs}
\Psi : K\Gamma \to \End_{K}(V^{\otimes r}); \qquad
\Phi:  K\Sym_r  \to \End_{K}(V^{\otimes r}) 
\end{equation}
induce inclusions
\begin{equation}\label{inclusions}
\Psi( K\Gamma ) \subseteq \End_{\Sym_r}(V^{\otimes r}); \qquad
\Phi( K\Sym_r ) \subseteq \End_{\Gamma}(V^{\otimes r}) 
\end{equation}
where $\End_{\Sym_r}(V^{\otimes r})$ (respectively,
$\End_{\Gamma}(V^{\otimes r})$) is defined to be the algebra of linear
operators on $V^{\otimes r}$ commuting with all operators in $\Phi(
\Sym_r )$ (respectively, $\Psi( \Gamma )$).  Equivalently, the
commutativity of the two actions says that the representations in
\eqref{eq:repdefs} induce algebra homomorphisms
\begin{equation}\label{induced}
\overline{\Psi}: K\Gamma  \to \End_{\Sym_r}(V^{\otimes r}); \qquad
\overline{\Phi}: K\Sym_r \to \End_{\Gamma}(V^{\otimes r}).
\end{equation}
The statement that has come to be known as ``Schur--Weyl duality'' is
the following.

\begin{thm}[Schur--Weyl duality] \label{thm:SWD}
  For any infinite field $K$, the inclusions in \eqref{inclusions} are
  actually equalities. Equivalently, the induced maps in
  \eqref{induced} are surjective.
\end{thm}

In case $K=\C$ this goes back to a classic paper of Schur
\cite{Schur}.\footnote{A proof of Schur--Weyl duality over $\C$ can be
  extracted from Weyl's book \cite{Weyl}. A detailed and accessible
  proof is written out in \cite[Theorem 3.3.8]{GoodmanWallach}.}  The main
purpose of this survey is to write out a complete proof of the theorem
for an arbitrary infinite field, \emph{assuming} the truth of the
  result in case $K=\C$. The strategy, suggested by S.~Koenig, is to
argue that the dimension of each of the four algebras in the
inclusions \eqref{inclusions} is independent of the characteristic of
the infinite field $K$. The claim for a general infinite field $K$
then follows immediately from the classical result over $\C$, by
dimension comparison.

We make no claim that this strategy is ``best'' in any sense; it is
merely one possible approach. For a completely different recent
approach, see \cite{KSX}.

\section{Surjectivity of  $\overline{\Psi}$}
\noindent
Let us first establish half of Theorem~\ref{thm:SWD}, namely the
surjectivity of the induced map $\overline{\Psi}:K\Gamma
\to \End_{\Sym_r}(V^{\otimes r})$ in \eqref{induced}.  For a very
direct (and shorter) approach to this result, see the argument on page
210 of \cite{CL}.  As already stated, the strategy followed here is to
argue that the algebras $\Psi(K\Gamma)$, $\End_{\Sym_r}(V^{\otimes
  r})$ have dimension (as vector spaces over $K$) which is independent
of the characteristic of the infinite field $K$.

We first establish that $\dim_K \Psi(K\Gamma)$ is independent of $K$
(so long as $K$ is infinite). For this we need a general principal,
which states that the ``envelope'' and ``coefficient space'' of a
representation are dual to one another.  To formulate the principle,
let $\Gamma$ be any semigroup and $K$ any field (not necessarily
infinite). Denote by $K^\Gamma$ the $K$-algebra of $K$-valued
functions on $\Gamma$, with the usual product and sum of elements
$f,f'$ of $K^\Gamma$ given by $(ff')(g) = f(g)f'(g)$,
$(f+f')(g)=f(g)+f'(g)$, for $g \in \Gamma$.

Given a representation $\tau: \Gamma \to \End_K(M)$ in a $K$-vector
space $M$, the {\em coefficient space} of the representation is by
definition the subspace $\cf_\Gamma M$ of $K^\Gamma$ spanned by the
coefficients $\{ r_{ab} \}$ of the representation. The coefficients
$r_{ab} \in K^\Gamma$ are determined relative to a choice of basis
$v_a$ ($a \in I$) for $M$ by the equations
\begin{equation}\label{gen:a}
\tau(g)\, v_b = \sum_{a \in I} r_{ab}(g)\, v_a
\end{equation}
for $g \in \Gamma$, $b \in I$. 

Let $K\Gamma$ be the semigroup algebra of $\Gamma$.  Elements of
$K\Gamma$ are sums of the form $\sum_{g\in \Gamma} a_g g$ ($a_g \in
K$) with finitely many $a_g \ne 0$. The group multiplication extends
by linearity to $K\Gamma$.  The given representation $\tau: \Gamma
\to \End_K(M)$ extends by linearity to an algebra homomorphism
$K\Gamma \to \End_K(M)$; by abuse of notation we denote this extended
map also by $\tau$.  The \emph{envelope}\footnote{This terminology is
  adapted from \cite{Weyl}, where Weyl writes about
  the ``enveloping algebra'' of a group representation as the algebra
  generated by the endomorphisms on the representing space coming from
  the action of all group elements. In modern terminology, this is
  just the image of the representation's linear extension to the group
  algebra.}  of the representation $\tau$ is by definition the
subalgebra $\tau(K\Gamma)$ of $\End_K(M)$.  The representation $\tau$
factors through its envelope; that is, we have a commutative diagram
\begin{equation}
\begin{gathered}
\xymatrix{
K\Gamma \ar[rr]^\tau \ar@{->>}[dr] && \End_K(M)\\
& \tau(K\Gamma) \ar@{^{(}->}[ur]
}
\end{gathered}
\end{equation}
in which the leftmost and rightmost diagonal arrows are a surjection
and injection, respectively. Taking linear duals, the above
commutative diagram induces another one
\begin{equation}
\begin{gathered}
\xymatrix{
(K\Gamma)^*   && \End_K(M)^* \ar[ll]_{\tau^*} \ar@{->>}[dl] \\
& {\tau(K\Gamma)}^* \ar@{_{(}->}[ul]
}
\end{gathered}
\end{equation}
in which the leftmost and rightmost diagonal arrows are now an
injection and surjection, respectively.  There is a natural
isomorphism of vector spaces $(K\Gamma)^* \simeq K^\Gamma$, given by
restricting a linear $K$-valued map on $K\Gamma$ to $\Gamma$; its
inverse is given by the process of linearly extending a $K$-valued map
on $\Gamma$ to $K\Gamma$.

\begin{lem} [{\cite[Lemma 1.2]{DD}}] \label{gen:lem}
The coefficient space $\cf_\Gamma(M)$ may be identified with the image
of $\tau^*$, so there is an isomorphism of vector spaces
$(\tau(K\Gamma))^* \simeq \cf_\Gamma M$. 
\end{lem}

\begin{proof}
Relative to the basis $v_a$ ($a \in I$) the
algebra $\End_K(M)$ has basis $e_{ab}$ ($a,b \in I$), where $e_{ab}$
is the linear endomorphism of $M$ taking $v_b$ to $v_a$ and taking all
other $v_c$, for $c \ne b$, to 0. In terms of this notation, equation
\eqref{gen:a} is equivalent with the equality
\begin{equation}
\tau(g) = \sum_{a,b \in I} r_{ab}(g) \,e_{ab}.
\end{equation}
Let $e'_{ab}$ be the basis of $\End_K(M)^*$ dual to the basis
$e_{ab}$, so that $e'_{ab}$ is the linear functional on $\End_K(M)$
taking the value 1 on $e_{ab}$ and taking the value 0 on all other
$e_{cd}$. Then one checks that $\tau^*$ carries $e'_{ab}$ onto
$r_{ab}$.  This proves that $\cf_\Gamma(M)$ may be identified with the
image of $\tau^*$, as desired.
\end{proof}

We apply the preceding lemma to the representation $M=V^{\otimes r}$
of $\Gamma = \GL(V)$, to conclude that $\dim_K \Psi(K\Gamma)$ is equal
to $\dim_K \cf_\Gamma (V^{\otimes r})$. Now the reader may easily
check that coefficient spaces are \emph{multiplicative}, i.e.,
$\cf_\Gamma(M \otimes N) = \cf_\Gamma (M) \cdot \cf_\Gamma(N)$. Here
the multiplication takes place in $K^\Gamma$. We will apply this fact
to compute the dimension of $\cf_\Gamma (V^{\otimes r}) =
(\cf_\Gamma(V))^r$.

From now on we choose (and fix) a basis $\{v_1, \dots, v_n\}$ of $V$
and identify $V$ with $K^n$ and $\Gamma$ with $\GL_n(K)$, by means of
the chosen basis. Then the action of $\Gamma$ on $V$ is by matrix
multiplication.

\begin{lem}\label{lem:AKnr}
  For $\Gamma =\GL_n(K)$ and $K$ any infinite field,
  $\cf_\Gamma(V^{\otimes r})$ is the vector space $A_K(n,r)$
  consisting of all homogeneous polynomial functions on $\Gamma$ of
  degree $r$. We have $\dim_K A_K(n,r) = \binom{n^2+r-1}{r} = \dim_K
  \Psi(K\Gamma)$.
\end{lem}

\begin{proof}
  Let $c_{ij} \in K^\Gamma$ be the function which maps a matrix $g \in
  \Gamma$ onto its $(i,j)$th matrix entry.  By definition, a function
  $f \in K^\Gamma$ is polynomial\footnote{The notion of ``polynomial''
    functions on general linear groups goes back (at least) to Schur's
    1901 dissertation.} if it belongs to the polynomial algebra $K[
  c_{ij}: 1\le i,j \le n ]$. The $c_{ij}$ are algebraically
  independent since $K$ is infinite. Note that the $c_{ij}$ are the
  coefficients of $\Gamma$ on $V$, i.e., $\cf_\Gamma V = \sum_{1\le i,j
    \le n} Kc_{ij}$.

  An element $f \in K[ c_{ij}: 1\le i,j \le n ]$ is homogeneous of
  degree $r$ if $f(a g) = a^r f(g)$ for all $a\in K$ and all $g\in
  \Gamma$. Here we define $ag$ to be the matrix obtained from $g$ by
  multiplying each entry by the scalar $a$.

  Now from the equality $\cf_\Gamma V = \sum_{1\le i,j \le n} Kc_{ij}$
  and the multiplicativity of coefficient spaces, it follows that
  $\cf_\Gamma(V^{\otimes r})$ is the vector space $A_K(n,r)$
  consisting of all homogeneous polynomial functions on $\Gamma$ of
  degree $r$. The equality $\dim_K A_K(n,r) = \binom{n^2+r-1}{r}$, now
  follows by an easy dimension count (or one can look at
  \cite[\S2.1]{Green}), and this is the same as $\dim_K \Psi(K\Gamma)$
  by Lemma \ref{gen:lem}.
\end{proof}

The preceding lemma establishes the fact that $\dim_K
\cf_\Gamma(V^{\otimes r})$ is independent of the characteristic of $K$
(so long as $K$ is infinite). So we turn now to the task of
establishing a similar independence statement for
$\dim_K \End_{\Sym_r}(V^{\otimes r})$.
 
Let us restrict the action of $\Gamma$ to the ``maximal torus'' $T
\subset \Gamma$ given by all diagonal matrices in $\Gamma =
\GL_n(K)$. The abelian group $T$ is isomorphic to the direct product
$(K^\times) \times \cdots \times (K^\times)$ of $n$ copies of the
multiplicative group $K^\times$ of the field $K$, so its irreducible
representations are one-dimensional, given on a basis element $z$ by
the rule $\mathrm{diag}(a_1, \dots, a_n) (z) = a_1^{\lambda_1} \cdots
a_n^{\lambda_n} z$, for various $\lambda_i \in \N$.  For convenience
of notation, write $t = \mathrm{diag}(a_1, \dots, a_n)$, $\lambda =
(\lambda_1, \dots, \lambda_n)$, and $t^\lambda = a_1^{\lambda_1}
\cdots a_n^{\lambda_n}$. Now $T$ acts semisimply on $V^{\otimes r}$,
and we have a ``weight space decomposition''
\begin{equation} \label{eq:wtspace} \textstyle
  V^{\otimes r} = \bigoplus_{\lambda\in \N^n} (V^{\otimes r})_\lambda
\end{equation}
where $(V^{\otimes r})_\lambda = \{ m \in V^{\otimes r}: t m =
t^\lambda\, m, \text{ for all } t \in T \}$.

Since the action of $T$ on $V^{\otimes r}$ commutes with the place
permutation action of $\Sym_r$, it follows that each weight space
$(V^{\otimes r})_\lambda$ is a $K\Sym_r$-module. It is easy to write
out a basis for $(V^{\otimes r})_\lambda$ in terms of the given basis
$\{v_1, \dots, v_n\}$ of $V$. Clearly $V^{\otimes r}$ has a basis
consisting of simple tensors of the form $v_{i_1} \otimes \cdots
\otimes v_{i_r}$ for various multi-indices $(i_1, \dots, i_r)$
satisfying the condition $i_j \in \{1, \dots, n\}$ for each $1 \le j
\le r$. Each simple tensor $v_{i_1} \otimes \cdots \otimes v_{i_r}$
has weight $\lambda = (\lambda_1, \dots, \lambda_n)$ where $\lambda_i$
counts the number of indices $j$ such that $i_j = i$. Thus it follows
that $\sum_i \lambda_i = r$. Let us write $\Lambda(n,r)$ for the set
of all $\lambda \in \N^n$ such that $\sum_i \lambda_i = r$.  Then each
summand $(V^{\otimes r})_\lambda$ is zero unless $\lambda \in
\Lambda(n,r)$, so we may replace $\N^n$ by $\Lambda(n ,r)$ in the
decomposition \eqref{eq:wtspace}.

From the above it follows that a basis of $(V^{\otimes r})_\lambda$,
for any $\lambda \in \Lambda(n,r)$, is given by the set of all
$v_{i_1} \otimes \cdots \otimes v_{i_r}$ of weight $\lambda$.

As a $K\Sym_r$-module, the weight space $(V^{\otimes r})_\lambda$ may
be identified with a ``permutation'' module $M^\lambda$.  Typically,
$M^\lambda$ is defined as the induced module $\mathbf{1}
\otimes_{(K\Sym_\lambda)} (K\Sym_r)$, where by $\mathbf{1}$ we mean
the one dimensional module $K$ with trivial action, and where
$\Sym_\lambda$ is the Young subgroup
\[
\Sym_{ \{1, \dots, \lambda_1\} } \times \Sym_{
\{\lambda_1+1, \dots, \lambda_1+\lambda_2\} } \times \cdots
\times \Sym_{
\{\lambda_{n-1}+1, \dots, \lambda_{n-1}+\lambda_n\} }
\]
of $\Sym_r$ determined by $\lambda = (\lambda_1, \dots, \lambda_n)$.
By \cite[\S12D]{CR} this has a basis (over $K$) indexed by any set of
right\footnote{Reference \cite{CR} works with left modules instead of
  right ones, so for our purposes left and right need to be
  interchanged there.} coset representatives of $\Sym_\lambda$ in
$\Sym_r$.

\begin{lem}\label{lem:EndSymr}
  For any field $K$, $\dim_K \End_{\Sym_r} (V^{\otimes r})$ is
  independent of $K$. 
\end{lem}

\begin{proof}
  From the decomposition \eqref{eq:wtspace} it follows that we have a
  direct sum decomposition of $\End_{\Sym_r}(V^{\otimes r}) =
  \Hom_{\Sym_r}(V^{\otimes r}, V^{\otimes r})$ of the form
  \begin{equation*} \textstyle
  \End_{\Sym_r}(V^{\otimes r}) = \bigoplus_{\lambda, \mu \in
    \Lambda(n,r)} \Hom_{\Sym_r}((V^{\otimes r})_\lambda, (V^{\otimes
    r})_\mu).
  \end{equation*}
  By Lemma \ref{lem:Perm}(b) in the next section, we may identify
  \[
  \Hom_{\Sym_r}((V^{\otimes r})_\lambda, (V^{\otimes r})_\mu) \simeq
  \Hom_{\Sym_r}(M^\lambda, M^\mu)
  \] 
  for any $\lambda, \mu \in \Lambda(n,r)$. By Mackey's theorem (see
  \cite[\S44]{CR} or combine \cite[Proposition 22]{Serre} with
  Frobenius reciprocity), it follows that $\dim_K
  \Hom_{\Sym_r}(M^\lambda, M^\mu)$ is equal to the number of
  $(\Sym_\lambda, \Sym_\mu)$-double cosets in $\Sym_r$, which is
  independent of $K$. This proves the claim.  Alternatively, one can
  avoid the Mackey theorem by applying James \cite[Theorem
    13.19]{James} directly (see also \cite[Proposition 3.5]{DEH}).
\end{proof}

Now we can obtain the main result of this section, which proves half
of Schur--Weyl duality in positive characteristic. We remind the
reader that the validity of Theorem \ref{thm:SWD} for $K=\C$ is
assumed, so in particular $\Psi(\C\Gamma) =
\End_{\Sym_r}((\C^n)^{\otimes r})$.

\begin{prop}\label{prop:part1}
  For any infinite field $K$, the image $\Psi(K\Gamma)$ of the
  representation $\Psi$ is equal to the centralizer algebra
  $\End_{\Sym_r}(V^{\otimes r})$, so the map $\overline{\Psi}$ in
  \eqref{induced} is surjective.
\end{prop}

\begin{proof}
  By Lemmas \ref{lem:AKnr} and \ref{lem:EndSymr} we have equalities
\begin{gather*}
\dim_K \Psi(K\Gamma) = \dim_\C \Psi(\C\Gamma), \\
\dim_K \End_{\Sym_r}((K^n)^{\otimes r}) =
\dim_\C \End_{\Sym_r}((\C^n)^{\otimes r})
\end{gather*}
for any infinite field $K$.  Since $\Psi(\C\Gamma) =
\End_{\Sym_r}((\C^n)^{\otimes r})$ it follows that $\dim_K
\Psi(K\Gamma) = \dim_K \End_{\Sym_r}((K^n)^{\otimes r})$ for any
infinite field $K$, and thus by comparison of dimensions the first
inclusion in \eqref{inclusions} must be an equality.  Equivalently,
the map $\overline{\Psi}$ in \eqref{induced} is surjective.
\end{proof}

\section{Surjectivity of $\overline{\Phi}$}\noindent
It remains to establish the surjectivity of the induced map $\overline{\Phi}$ in
\eqref{induced}. This surjectivity was first established in positive
characteristic in \cite[Theorem 4.1]{De-Pr}.\footnote{The statement of
  Theorem 4.1 in \cite{De-Pr} is actually much more general.}  We will
outline an alternative proof here, following our avowed strategy of
showing that the dimensions of $\Phi(K\Sym_r)$,
$\End_{\Gamma}(V^{\otimes r})$ are independent of the characteristic
of the infinite field $K$.

In order to establish the independence statement for $\Phi(K\Sym_r)$
we apply results of Murphy and H\"{a}rterich in order to compute the
annihilator of the action of $\Sym_r$ on $V^{\otimes r}$. Note that
Murphy and H\"{a}rterich worked with the Iwahori--Hecke algebra (with
parameter $q$) in type $A$, so one needs to take $q=1$ in their
formulas in order to get corresponding results for the group algebra
$K\Sym_r$. The results of Murphy and H\"{a}rterich hold over an
arbitrary commutative integral domain, so $K$ does not need to be an
infinite field in this part. So we assume from now on, until the
paragraph after Corollary \ref{cor:Phi-invariance}, that $K$ is a
commutative integral domain.

Let $\lambda$ be a composition of $r$.  We regard $\lambda$ as an
infinite sequence $(\lambda_1, \lambda_2, \dots)$ of nonnegative
integers such that $\sum \lambda_i = r$. The individual $\lambda_i$
are the parts of $\lambda$, and the largest index $\ell$ such that
$\lambda_\ell = 0$ and $\lambda_j = 0$ for all $j > \ell$ is the
length, or number of parts, of $\lambda$. Any composition $\lambda$
may be sorted into a partition $\lambda^+$, in which the parts are
non-strictly decreasing. When writing compositions or partitions,
trailing zero parts are usually omitted. If $\lambda$ is a partition,
we generally write $\lambda'$ for the \emph{transposed} (or
\emph{conjugate}) partition, corresponding to writing the rows of the
Young diagram as columns.

Given a composition $\lambda=(\lambda_1, \dots, \lambda_\ell)$ of $r$,
a Young diagram of shape $\lambda$ is an arrangement of boxes into
rows with $\lambda_i$ boxes in the $i$th row.  A $\lambda$-tableau $T$
is a numbering of the boxes in the Young diagram of shape $\lambda$ by
the numbers $1, \dots, r$ so that each number appears just once. In
other words, it is a bijection between the boxes in the Young diagram
and the set $\{1, \dots, r\}$.  Such a $T$ is {\em row standard} if
the numbers in each row are increasing when read from left to right,
and {\em standard} if row standard and the numbers in each column are
increasing when read from top to bottom.

The group $\Sym_r$ acts naturally on tableaux, on the right, by
permuting the entries.  Given a tableau $T$, we define the {\em row
  stabilizer} of $T$ to be the subgroup $R(T)$ of $\Sym_r$ consisting
of those permutations that permute entries in each row of $T$ amongst
themselves, similarly the {\em column stabilizer} is the subgroup
$C(T)$ consisting of those permutations that permute entries in each
column of $T$ amongst themselves.

Let $\lambda$ be a composition of $r$. Let $T^\lambda$ be the
$\lambda$-tableau in which the numbers $1, \dots, r$ have been
inserted in the boxes in order from left to right along rows, read
from top to bottom. Set $\Sym_\lambda=R(T^\lambda)$. This is the same
as the Young subgroup
\[
\Sym_{ \{1, \dots, \lambda_1\} } \times \Sym_{
\{\lambda_1+1, \dots, \lambda_1+\lambda_2\} } \times \cdots
\]
of $\Sym_r$ defined by the composition $\lambda$.  Given a row
standard $\lambda$-tableau $T$, we define $d(T)$ to be the unique
element of $\Sym_r$ such that $T = T^\lambda d(T)$.  Given any pair
$S,T$ of row standard $\lambda$-tableaux, following Murphy
\cite{Murphy} we set
\begin{equation}
 x_{ST} = d(S)^{-1} x_\lambda d(T); \quad 
 y_{ST} = d(S)^{-1} y_\lambda d(T).
\end{equation}
where $x_\lambda = \sum_{w \in \Sym_\lambda} w$ and $y_\lambda =
\sum_{w \in \Sym_\lambda} (\sgn w)\,w$. 

\begin{thm}[Murphy]
  Let $K$ be a commutative integral domain. Each of the sets
  $\{x_{ST}\}$ and $\{y_{ST}\}$, as $(S,T)$ ranges over the set of all
  ordered pairs of standard $\lambda$-tableaux for all partitions
  $\lambda$ of $r$, is a $K$-basis of the group algebra $A=K\Sym_r$.
\end{thm}

Note that $x_{ST}$ and $y_{ST}$ are interchanged by the $K$-linear
ring involution of $K\Sym_r$ which sends $w$ to $(\sgn w) w$, for $w
\in \Sym_r$. This gives a trivial way of converting results about one
basis into results about the other.

We will need several equivalent descriptions of the permutation
modules $M^\lambda$, which we now formulate. Let $\lambda$ be a
composition of $r$. Recall that $M^\lambda = \mathbf{1}
\otimes_{(K\Sym_\lambda)} (K\Sym_r)$, where $\mathbf{1}$ is the one
dimensional module $K$ with trivial action. In \cite[Definition
  4.1]{James}, an alternative combinatorial description of $M^\lambda$
is given in terms of ``tabloids'' (certain equivalence classes of
tableaux), and in \cite[(1.3)]{DJ} the authors write out an explicit
isomorphism between these two descriptions. The following gives two
additional descriptions of $M^\lambda$, the second of which was used
already in the previous section.

\begin{lem} \label{lem:Perm}
  For any composition $\lambda$ of $r$, the permutation module
  $M^\lambda$ is isomorphic (as a right $K\Sym_r$-module) with either
  of 

  (a) the right ideal $x_\lambda (K\Sym_r)$ of $K\Sym_r$;

  (b) the weight space $(V^{\otimes r})_\lambda$ in $V^{\otimes r}$,
  where $V$ is free over $K$ of rank at least as large as the number
  of parts of $\lambda$.
\end{lem}

\begin{proof} Let $\mathcal{D}_\lambda = \{ d(T) \}$ as $T$ varies
  over the set of row standard tableaux of shape $\lambda$. This is a
  set of right coset representatives of $\Sym_\lambda$ in
  $\Sym_r$. The map $d \to x_\lambda d$ gives the isomorphism (a), in
  light of Lemma 3.2(i) of \cite{DJ}. The isomorphism (b) works as
  follows. Given $d \in \mathcal{D}_\lambda$, write $d=d(T)$ for some
  (unique) row standard tableau $T$ of shape $\lambda$. Use $T$ to
  construct a simple tensor $v_{i_1} \otimes \cdots \otimes v_{i_r}$
  of weight $\lambda$, by letting $i_j$ be the (unique) row number in
  $T$ in which $j$ is found. This map is well defined, and is a
  bijection since there is an obvious inverse map.
\end{proof}

We recall that compositions are partially ordered by \emph{dominance},
defined as follows. Given two compositions $\lambda, \mu$ of $r$,
write $\lambda \dom \mu$ ($\lambda$ dominates $\mu$) if $\sum_{i\le j}
\lambda_i \ge \sum_{i\le j} \mu_i$ for all $j$. One writes $\lambda
\sdom \mu$ ($\lambda$ strictly dominates $\mu$) if $\lambda \dom \mu$
and the inequality $\sum_{i\le j} \lambda_i \ge \sum_{i\le j} \mu_i$
is strict for at least one $j$.

The dominance order on compositions extends to the set of row standard
tableaux, as follows. Let $T$ be a row standard $\lambda$-tableau,
where $\lambda$ is a composition of $r$. For any $s<r$ denote by
$T_{\downarrow s}$ the row standard tableau that results from throwing
away all boxes of $T$ containing a number bigger than $s$. Let
$[T_{\downarrow s}]$ be the corresponding composition of $s$ (the
composition defining the shape of $T_{\downarrow s}$). Given row
standard tableaux $S,T$ with the same number $r$ of boxes, define
\begin{equation}\label{eq:domtab}
\begin{gathered}
 \text{$S \dom T$ if for each $s \le r$, $[S_{\downarrow s}] \dom
 [T_{\downarrow s}]$};\\
 \text{$S \sdom T$ if for each $s \le r$, $[S_{\downarrow s}] \sdom
 [T_{\downarrow s}]$}.
\end{gathered}
\end{equation}
Note that if $S,T$ are standard tableaux, respectively of shape
$\lambda,\mu$ where $\lambda$ and $\mu$ are partitions of $r$, then $S
\dom T$ if and only if $T' \dom S'$. Here $T'$ denotes the transposed
tableau of $T$, obtained from $T$ by writing its rows as columns.

Let $*$ be the $K$-linear anti-involution on $A=K\Sym_r$ given by 
\[ \textstyle
(\sum_{w\in \Sym_r} b_w w)^* \to \sum_{w\in \Sym_r} b_w w^{-1}
\]
for any $b_w \in K$.  An easy calculation with the definitions shows
that
\begin{equation} \label{eq:star}
x_{ST}^{*} = x_{TS}^{\ }; \quad  y_{ST}^{*} = y_{TS}^{\ }
\end{equation}
for any pair $S,T$ of row standard $\lambda$-tableaux.

We write $c \in \{x,y\}$ in order to describe the cell structure of
$A=K\Sym_r$ relative to both bases simultaneously.  

\begin{thm} [Murphy, {\cite[Theorem 4.18]{Murphy}}]
  Let $c \in \{x,y\}$. Let $\lambda$ be a partition of $r$. The
  $K$-module $A[\dom \lambda] = \sum K c_{ST}$, the sum taken over all
  pairs $(S,T)$ of standard $\mu$-tableaux such that $\mu \dom
  \lambda$, is a two-sided ideal of $A$, as is $A[\sdom \lambda] =
  \sum K c_{ST}$, the sum taken over all pairs $(S,T)$ of standard
  $\mu$-tableaux such that $\mu \sdom \lambda$. For any $a \in A$ and
  any pair $(S,T)$ of $\lambda$-tableaux, we have
\begin{equation} \label{eq:cellstructure} \textstyle
  c_{ST}\, a = \sum_{U} r_a(T,U)\, c_{SU} \mod A[\sdom \lambda]
\end{equation}
where $r_a(T,U) \in K$ is independent of $S$, and in the sum $U$
varies over the set of standard $\lambda$-tableaux. 
\end{thm}

In the language of cellular algebras, introduced by Graham and Lehrer
\cite{GL}, for $c \in \{x,y\}$ the basis $\{c_{ST}\}$ is a cellular
basis of $A$. Note that by applying the anti-involution $*$ to
\eqref{eq:cellstructure} we obtain by \eqref{eq:star} the equivalent
condition
\begin{equation} \label{eq:dualcellstructure} \textstyle
  a^*\, c_{TS} = \sum_{U} r_a(T,U)\, c_{US} \mod A[\sdom \lambda]
\end{equation}
for any $a \in A$ and any pair $(S,T)$ of $\lambda$-tableaux.

Now fix $n$ and $r$, and let $P$ be the set of partitions $\lambda$ of
$r$ such that $\lambda_1 > n$.  Note that $P$ is empty if $n \ge r$.
Set $A[P]= \sum K y_{ST}$, where the sum is taken over the set of
pairs $(S,T)$ of standard tableaux of shape $\lambda$, for all
$\lambda \in P$.  It follows from \eqref{eq:cellstructure},
\eqref{eq:dualcellstructure} that $A[P]$ is a two-sided ideal of $A$
because $P$ satisfies the property: $\lambda \in P$, $\mu \dom
\lambda$ $\Longrightarrow$ $\mu \in P$ for any partition $\mu$ of $r$.
Note that $A[P]$ is the zero ideal if $n \ge r$.

\begin{lem} \label{lem:contain}
The kernel of $\Phi$ contains $A[P]$.
\end{lem}

\begin{proof}
  If $n \ge r$ then $P$ is empty and there is nothing to prove, so we
  may assume that $n < r$.

  We first observe that $y_\lambda$ acts as zero on any simple tensor
  $v_{i_1} \otimes \cdots \otimes v_{i_r} \in V^{\otimes r}$, for any
  $\lambda \in P$. This is because any such tensor has at most $n$
  distinct tensor factors, and thus is annihilated by the alternating
  sum $\alpha=\sum_{w \in \Sym_{\{1, \dots, \lambda_1\}}} (\sgn
  w)w$. (Recall that $\lambda_1 > n$.) The alternating sum $\alpha$ is
  a factor of $y_\lambda$, i.e., we have $y_\lambda = \alpha \beta$
  for some $\beta \in K\Sym_r$, so $y_\lambda$ acts as zero as
  well. Since $V^{\otimes r}$ is spanned by such simple tensors, it
  follows that $y_\lambda$ acts as zero on $V^{\otimes r}$.

  It follows immediately that every $y_{ST} = d(S)^{-1} y_\lambda
  d(T)$, for $\lambda \in P$, acts as zero on $V^{\otimes r}$, for any
  $\lambda$-tableaux $S,T$, since $d(S)^{-1}$ simply permutes the
  entries in the tensor, and then $y_\lambda$ annihilates it. Since
  $A[P]$ is spanned by such $y_{ST}$, it follows that $A[P]$ is
  contained in the kernel of $\Phi$.
\end{proof}

We will use a lemma of Murphy to establish the opposite inclusion. Let
$(S,T)$ be a pair of $\lambda$-tableaux, where $\lambda$ is a
composition of $r$. The pair is row standard if both $S,T$ are row
standard; similarly the pair is standard if both $S,T$ are standard.
The dominance order on tableaux defined in \eqref{eq:domtab} extends
naturally to pairs of tableaux, by defining:
\begin{equation}
    (S,T) \dom (U,V) \text{ if } S\dom U \text{ and } T\dom V.
\end{equation}
For $a,b \in A$ let $(a,b)$ denote the coefficient of 1 in the
expression $ab^{*} = \sum_{w\in \Sym_r} c_w\, w$, where $c_w \in
K$. Then $(\ ,\ )$ is a non-degenerate symmetric bilinear form on
$A=K\Sym_r$. It is straightforward to check that this bilinear form
satisfies the properties
\begin{equation} \label{eqn:bilform}
  (a, bd) = (ad^{*}, b); \quad (a, db) = (d^{*}a, b)
\end{equation}
for any $a,b,d \in A$.

\begin{lem}[{Murphy, \cite[Lemma 4.16]{Murphy2}}] \label{lem:keylem} 
Let $(S,T)$ be a row standard pair of $\mu$-tableaux and $(U,V)$ a
standard pair of $\lambda$-tableaux, where $\mu$ is a given composition
of $r$ and $\lambda$ a partition of $r$. Then:

(a) $(x_{ST}, y_{U'V'}) = 0$ unless $(U,V) \dom (S,T)$;

(b) $(x_{UV}, y_{U'V'}) = \pm 1$

\noindent
where $T'$ denotes the transpose of a tableau $T$. 
\end{lem}

\noindent
This is used in proving the following result, which in particular
shows that the rank (over $K$) of the annihilator of the symmetric
group action on $V^{\otimes r}$ is independent of the characteristic
of $K$.

\begin{prop}[{H\"{a}rterich, \cite[Lemma~3]{Harterich}}] \label{thm:main}
The kernel of $\Phi$, i.e., the annihilator $\ann^{\ }_{K\Sym_r} V^{\otimes
r}$, is the cell ideal $A[P]$.
\end{prop}

\begin{proof}
  By Lemma \ref{lem:contain}, the kernel of $\Phi$ contains $A[P]$, so
  we only need to prove the reverse containment.  Let
\[ \textstyle
a = \sum_{(S,T)} a_{ST}y_{ST} \in \ker
\Phi
\]
where $a_{ST} \in K$, and the sum over all pairs $(S,T)$ of standard
tableaux of shape $\lambda$, where $\lambda$ is a partition of $r$.
It suffices to prove: ($*$) $a_{ST} = 0$ for all pairs $(S,T)$ of
standard tableaux of shape $\mu \in P^c$, where $P^c$ is the
complement of $P$ in the set of all partitions of $r$.

We note that $P^c$ is the set of conjugates $\lambda'$ of partitions
$\lambda$ in $\Lambda(n,r)$. Write $\Lambda^+(n,r)$ for the set of
partitions in $\Lambda(n,r)$; this is the set of partitions of $r$
into not more than $n$ parts.

We proceed by contradiction. Suppose $(*$) is not true. Since by Lemma
\ref{lem:contain} we have $\sum_{\shape(S,T)\in P} a_{ST}y_{ST} \in
\ker(\Phi)$, it follows that
\[ \textstyle
  b = \sum_{\shape(S,T)\in P^c} a_{ST}y_{ST} 
\]
is also in the kernel of $\Phi$; i.e., the element $b$ annihilates
$V^{\otimes r}$.  Under the assumption we have $b\ne 0$.  Let
$(S_0,T_0)$ be a minimal pair (with respect to $\dom$) with
$\shape(S_0,T_0) \in P^c$ such that $a_{S_0T_0} \ne 0$. So $a_{ST} =
0$ for all pairs $(S,T)$ with $(S_0, T_0) \sdom (S,T)$. Let
$\lambda_0$ be the shape of $T'_0$ (same as shape of $S'_0$). Then
$\lambda_0 \in \Lambda^+(n,r)$, and we have
\begin{align*}
  ( x_{\lambda_0 S'_0}\, b, d(T'_0) ) &= ( x_{\lambda_0 S'_0}\,
     \sum  a_{ST}y_{ST} , d(T'_0) ) \\
  &= \sum a_{ST} ( d(T'_0)^{-1} x_{\lambda_0 S'_0}, y_{ST}^* ) \\
  &= \sum a_{ST} ( x_{T'_0 S'_0}, y_{TS} ) 
\end{align*}
where all sums are taken over the set of $(S,T)$ of shape some member
of $P^c$.  Here, we write $x_{\mu T}$ shorthand for $x_{T^\mu T}$,
where (as before) $T^\mu$ is the $\mu$-tableau in which the numbers
$1, \dots, r$ have been inserted in the boxes in order from left to
right along rows, read from top to bottom.

By Lemma \ref{lem:keylem}(a) all the terms in the last sum are zero
unless $(S_0,T_0) \dom (S,T)$, in other words $( x_{T'_0 S'_0}, y_{TS}
) = 0$ for all pairs $(S,T)$ which are strictly more dominant than
$(S_0,T_0)$. By assumption, $a_{ST}=0$ for all pairs $(S,T)$ strictly
less dominant than $(S_0,T_0)$. Thus, the above sum collapses to a
single term $a_{S_0T_0} ( x_{T'_0 S'_0}, y_{T_0S_0} )$, and by our
assumption and Lemma \ref{lem:keylem}(b) this is nonzero.

This proves that $x_{\lambda_0 S'_0}\, b \ne 0$. Thus $b$ does not
annihilate the permutation module $M^{\lambda_0} \simeq
x_{\lambda_0}A$.  Since $\lambda_0 \in \Lambda^+(n,r)$ as noted above,
and thus $M^{\lambda_0}$ is isomorphic to a direct summand of
$V^{\otimes r}$, we have arrived at a contradiction. This proves the
result.
\end{proof}

\begin{cor}\label{cor:Phi-invariance}
  For any commutative integral domain $K$, the $K$-module
  $\Phi(K\Sym_r)$ is free over $K$, of rank $r! - \sum_{\lambda \in P}
  N(\lambda)^2$, where $N(\lambda)$ is the number of standard tableaux
  of shape $\lambda$.  In particular, the $K$-rank of $\Phi(K\Sym_r)$
  is independent of $K$.
\end{cor}

\begin{proof}
  By the preceding proposition, $\Phi(K\Sym_r) \simeq A/A[P]$. This is
  free over $K$ because it is a submodule of the free $K$-module
  $\End_K(V^{\otimes r})$. By definition, $A[P]$ is free over $K$ of rank 
  $\sum_{\lambda \in P} N(\lambda)^2$, so the result follows.
\end{proof}

Now we return to the assumption that $K$ is an infinite field, and
consider why $\dim_K \End_{\Gamma}(V^{\otimes r})$ is independent of
$K$. This involves facts about the representation theory of algebraic
groups that are less elementary than facts used so far. We identify
the group $\Gamma = \GL_n(K)$, the group of $K$-rational points in the
algebraic group $\GL_n(\overline{K})$, where $\overline{K}$ is an
algebraic closure of $K$, with the group scheme $\mathbf{GL}_n$ over
$K$.

For $\Gamma = \GL_n(K)$ we let $T$ be the maximal torus consisting of
all diagonal elements of $\Gamma$. Regard an element $\lambda \in
\Z^n$ as a character on $T$ (via $\mathrm{diag}(a_1, \dots, a_n) \to
a_1^{\lambda_1} \cdots a_n^{\lambda_n}$ for $a_i \in
K^\times$). Consider the Borel subgroup $B$ consisting of the lower
triangular matrices in $\Gamma$, and let $\nabla(\lambda)$ be the
induced module (see \cite[Part I, \S3.3]{Jantzen}):
\[
\ind_B^\Gamma (K_\lambda) = \{f \in K[\Gamma]: f(gb) = b^{-1} f(g),
\text{ all } b\in B, g\in G \}
\]
for any $\lambda \in \Z^n$, where $K_\lambda$ is the one dimensional
$T$-module with character $\lambda$, regarded as a $B$-module by
making the unipotent radical of $B$ act trivially. 

The dual space $M^* = \Hom_K(M, K)$ of a given rational
$K\Gamma$-module $M$ is again a rational $K\Gamma$-module, in two
different ways:

(i) \quad $(g\cdot f)(m) = f(g^{-1} m)$;

(ii) \quad $(g\cdot f)(m) = f(g^{t} m)$ \qquad ($g^\transpose$ is the
matrix transpose of $g$)

\noindent
for $g\in \Gamma$, $f\in M^*$, $m \in M$. Denote the first dual by
$M^*$ and the second by $M^\transpose$.  Let $\Delta(\lambda) =
\nabla(\lambda)^\transpose$.  It is known that $\Delta(\lambda) \simeq
\nabla(-w_0 \lambda)^*$ where $w_0$ is the longest element in the Weyl
group $W$. The modules $\nabla(\lambda)$, $\Delta(\lambda)$ are known
as ``dual Weyl modules'' and ``Weyl modules'',
respectively.\footnote{Weyl and dual Weyl modules for $\GL_n(K)$ are
  studied in \cite[Chapters 4, 5]{Green}, where they are respectively
  denoted by $D_{\lambda, K}$ and $V_{\lambda, K}$.}  The most
important property these modules satisfy, for our purposes, is the
following
\begin{equation} \label{eq:Ext}
  \operatorname{Ext}^j_\Gamma (\Delta(\lambda), \nabla(\mu)) \simeq 
\begin{cases}
  K & \text{ if } j=0 \text{ and } \lambda = \mu \\ 
  0 & \text{ otherwise.}
\end{cases}
\end{equation}
This is a special case of \cite[Part II, Proposition 4.13]{Jantzen}. 

Say that a $\Gamma$-module $M$ has a $\nabla$-filtration
(respectively, $\Delta$-filtration) if it has an ascending chain of
submodules 
\[
0=M_0 \subseteq M_1 \subseteq \cdots \subseteq M_{t-1}
\subseteq M_t = M
\]
such that each successive quotient $M_i/M_{i-1}$ is isomorphic with
$\nabla(\lambda^i)$ (respectively, $\Delta(\lambda^i)$) for some
$\lambda^i \in \Z^n$. Another fact we need goes back to \cite[Theorem
B, page 164]{Wang}:
\begin{equation} \label{eq:WeylFiltration}
  \Delta(\lambda) \otimes \Delta(\mu) \text{ has a $\Delta$-filtration}
\end{equation}
for any $\lambda, \mu \in \Z^n$. (Note that this fundamental result
has been extended in \cite{Donkin:tensor}, which in turn was
extended in \cite{Mathieu}.)  From \eqref{eq:WeylFiltration} it
follows immediately by taking duals that
\begin{equation} \label{eq:goodFiltration}
  \nabla(\lambda) \otimes \nabla(\mu) \text{ has a $\nabla$-filtration}
\end{equation}
for any $\lambda, \mu \in \Z^n$. The following result, which says that
$V^{\otimes r}$ is a ``tilting'' module for $\Gamma$, is now easy to
prove.

\begin{lem}\label{lem:tensorFilt}
  $V^{\otimes r}$ has both $\nabla$- and $\Delta$-filtrations.
\end{lem}

\begin{proof}
One has $V = \nabla(\varepsilon_1) = \Delta(\varepsilon_1)$ where
$\varepsilon_1 = (1, 0, \dots, 0)$. The result then follows from
\eqref{eq:WeylFiltration} and \eqref{eq:goodFiltration} by induction
on $r$. 
\end{proof}

For the next argument we will need the notion of formal
characters. Any rational $K\Gamma$-module $M$ has a weight space
decomposition $M = \bigoplus_{\lambda \in \Z^n} M_\lambda$ where
\[
  M_\lambda = \{ m \in M: tm = t^\lambda m, \text{ for all } t\in T\}. 
\]
Here $t^\lambda = a_1^{\lambda_1} \cdots a_n^{\lambda_n}$ where $t =
\textrm{diag}(a_1, \dots a_n)$ as previously defined, just before
\eqref{eq:wtspace}. Set $X = \Z^n$ and let $\Z[X]$ be the free
$\Z$-module on $X$ with basis consisting of all symbols $e(\lambda)$
for $\lambda \in X$, with a multiplication given by $e(\lambda)e(\mu)
= e(\lambda+\mu)$, for $\lambda, \mu \in X$. If $M$ is finite
dimensional, the formal character $\ch M \in \Z[X]$ of $M$ is defined
by
\[ \textstyle
 \ch M = \sum_{\lambda \in X} (\dim_K M_\lambda)\, e(\lambda).
\]
The formal character of $\Delta(\lambda)$, which is the same as $\ch
\nabla(\lambda)$ since the maximal torus $T$ is fixed pointwise by the
matrix transpose, is given by Weyl's character formula \cite[Part II,
  Proposition 5.10]{Jantzen}.\footnote{The computation of the $\ch
  \Delta(\lambda)$ for $\GL_n(\C)$ goes back to Schur's 1901
  dissertation. Thus, these characters are sometimes called
  \emph{Schur functions.} See \cite{Macdonald} or \cite[Chapter
    7]{Stanley} for exhaustive accounts of their many properties.}

\begin{prop} \label{prop:EndGamma} For any infinite field $K$,
  $\dim_K \End_{\Gamma}(V^{\otimes r})$ is independent of $K$.
\end{prop}

\begin{proof}
  Let $0 = N_0 \subseteq N_1 \subseteq \cdots \subseteq N_{s-1}
  \subseteq N_s = V^{\otimes r}$ be a $\nabla$-filtration and $0 = M_0
  \subseteq M_1 \subseteq \cdots \subseteq M_{t-1} \subseteq M_t =
  V^{\otimes r}$ a $\Delta$-filtration. Write $(V^{\otimes r}:
  \nabla(\lambda))$ for the number of successive subquotients
  $N_i/N_{i-1}$ which are isomorphic to $\nabla(\lambda)$, and
  similarly write $(V^{\otimes r}: \Delta(\lambda))$ for the number of
  successive subquotients $M_i/M_{i-1}$ which are isomorphic to
  $\Delta(\lambda)$. Since characters are additive on short exact
  sequences, we have
  \[
   \ch V^{\otimes r} = \sum_{\lambda \in X} (V^{\otimes r}:
   \nabla(\lambda)) \ch \nabla(\lambda) = \sum_{\lambda \in X}
   (V^{\otimes r}: \Delta(\lambda)) \ch \Delta(\lambda).
  \]
  Since $V^{\otimes r}$ is self-dual (under the transpose dual) we may
  choose the filtration $(N_*)$ to be dual to the filtration
  $(M_*)$. It follows that $s=t$ and $(V^{\otimes r}: \nabla(\lambda))
  = (V^{\otimes r}: \Delta(\lambda))$ for all $\lambda$.

  Now one applies \eqref{eq:Ext} and a double induction through the
  filtrations. The argument is standard homological algebra, safely
  left at this point as an exercise for the reader. At the end one
  finds that
  \[ \textstyle 
  \dim_K \End_{\Gamma}(V^{\otimes r}) = \sum_{\lambda \in \Z^n}
  (V^{\otimes r}: \nabla(\lambda))^2
  \]
  where the number of nonzero terms in the sum is finite.  The result
  follows.
\end{proof}

Now we are ready to prove the second half of Schur--Weyl duality in
positive characteristic. We remind the reader that we assume the
validity of Theorem \ref{thm:SWD} in case $K=\C$.

\begin{prop} \label{prop:part2}
  For any infinite field $K$, the image $\Phi(K\Sym_r)$ of the
  representation $\Phi$ is equal to the centralizer algebra
  $\End_{\Gamma}(V^{\otimes r})$, so the map $\overline{\Phi}$ in
  \eqref{induced} is surjective.
\end{prop}

\begin{proof} The argument is essentially the same as the proof of
  Proposition \ref{prop:part1}.  By Corollary \ref{cor:Phi-invariance}
  and Proposition \ref{prop:EndGamma} we have equalities
\begin{gather*}
\dim_K \Phi(K\Sym_r) = \dim_\C \Phi(\C\Sym_r), \\
\dim_K \End_{\GL_n(K)}((K^n)^{\otimes r}) =
\dim_\C \End_{\GL(_n(\C)}((\C^n)^{\otimes r})
\end{gather*}
for any infinite field $K$.  Since $\Phi(\C\Sym_r) =
\End_{\GL_n(\C)}((\C^n)^{\otimes r})$ it follows that $\dim_K
\Phi(K\Sym_r) = \dim_K \End_{\GL_n(K)}((K^n)^{\otimes r})$ for any
infinite field $K$, and thus by comparison of dimensions the second
inclusion in \eqref{inclusions} must be an equality.  Equivalently,
the map $\overline{\Phi}$ in \eqref{induced} is surjective.
\end{proof}

By putting together Propositions \ref{prop:part1} and \ref{prop:part2}
we have now established Theorem \ref{thm:SWD} in positive
characteristic, assuming its validity for $K=\C$.

\begin{rmk}
(a) Let $K$ be an arbitrary infinite field. Lemma \ref{lem:AKnr} gives
  the equality $\dim_K(K\Gamma) = \binom{n^2+r-1}{r}$, and the proof
  of Lemma \ref{lem:EndSymr} in light of \cite[Theorem 13.19]{James}
  gives the equality $\dim_K \End_{\Sym_r}(V^{\otimes r}) =
\sum_{\lambda, \mu \in\Lambda(n,r)} N(\lambda^+, \mu^+)$, where
$N(\lambda^+, \mu^+)$ counts the number of ``semistandard'' tableaux
of shape $\lambda^+$ and weight $\mu^+$. Corollary
\ref{cor:Phi-invariance} says that $\dim_K \Phi(K\Sym_r) = r! -
\sum_{\lambda \in P} N(\lambda)^2$, where $N(\lambda)$ is the number
of standard tableaux of shape $\lambda$, and the proof of Proposition
\ref{prop:EndGamma} shows that $\dim_K \End_\Gamma(V^{\otimes r}) =
\sum_{\lambda \in \Lambda^+(n,r)} (V^{\otimes r}: \nabla(\lambda))^2$.
Thus, in order to obtain a proof of Theorem \ref{thm:SWD} in full
generality (without assuming its validity for $K=\C$) from the methods
of this paper, one only needs to demonstrate the combinatorial
identities
\begin{gather}
  \binom{n^2+r-1}{r} = \sum_{\lambda, \mu \in\Lambda(n,r)}
  N(\lambda^+, \mu^+);\\ r! - \sum_{\lambda \in P} N(\lambda)^2 =
  \sum_{\lambda \in \Lambda^+(n,r)} (V^{\otimes r}: \nabla(\lambda))^2 .
\end{gather}
The author has not attempted to construct a combinatorial proof of
these identities. If one assumes the validity of Theorem \ref{thm:SWD}
in the case $K=\C$, then these identities follow from the results in
this paper. Alternatively, if one can find an independent proof of the
identities, then one would have a new proof of Theorem \ref{thm:SWD}
in full generality, including the case $K=\C$.

(b) There is a variant of Theorem \ref{thm:SWD} worth noting. One may
twist the action of $\Sym_r$ on $V^{\otimes r}$ by letting $w\in
\Sym_r$ act as $(\sgn w)w$ (so $\Sym_r$ acts by ``signed'' place
permutations). This action also commutes with the action of $\Gamma =
\GL(V)$, and Theorem \ref{thm:SWD} holds for this action as well. This
may be proved the same way. In the course of carrying out the
argument, one needs to replace permutation modules by ``signed''
permutation modules, and interchange the role of Murphy's two bases
$\{ x_{ST} \}$, $\{ y_{ST} \}$.

(c) There is also a $q$-analogue of Theorem \ref{thm:SWD}, in which
one replaces $\GL_n(K)$ by the quantized enveloping algebra
corresponding to the Lie algebra $\gl_n$, and replaces $K\Sym_r$ by
the Iwahori--Hecke algebra $\mathbf{H}(\Sym_r)$. The generic case ($q$
not a root of unity) of this theorem was first observed in Jimbo
\cite{Jimbo}, and the root of unity case was treated in Du, Parshall,
and Scott \cite{DPS}. Alternatively, one may derive the result in the
root of unity case from Jimbo's generic version, using arguments along
the lines of those sketched here.
\end{rmk}

\bibliographystyle{amsalpha}

\end{document}